\documentclass[a4paper,11pt]{amsart}

\usepackage{amsthm,amsmath,amssymb,amsfonts}
\usepackage[utf8]{inputenc}
\usepackage{verbatim}
\usepackage{marvosym}
\usepackage{stackrel}
\usepackage{graphicx}
\usepackage[svgnames]{xcolor}
\definecolor{blue(munsell)}{rgb}{0.0, 0.5, 0.69}
\usepackage{lipsum}
\usepackage{svg}
\usepackage{caption}
\usepackage{enumitem}
\usepackage{epigraph}
\usepackage{ebproof}
\usepackage[cal=boondoxo]{mathalfa}
\usepackage[all,cmtip]{xy}
\usepackage{mathtools}
\usepackage{scalerel}
\usepackage{color}
\usepackage{mathrsfs}
\usepackage{subfig}
\usepackage{framed}
\usepackage[mathscr]{eucal}
\usepackage{hyperref}
\hypersetup{hypertexnames = false, bookmarksdepth = 2, bookmarksopen = true, colorlinks, linkcolor = black, citecolor = blue(munsell), urlcolor = blue(munsell), pdfstartview={XYZ null null 1}}
\usepackage{cleveref}
\usepackage{tikz}
\usepackage{quiver}

\usepackage{bbold}
\usepackage{breqn}
\usepackage{booktabs}

\usepackage{graphicx}

\usepackage[colorinlistoftodos]{todonotes}
\usepackage{xparse}
\usepackage{scalerel,stackengine}
\stackMath
\newcommand\reallywidehat[1]{%
	\savestack{\tmpbox}{\stretchto{%
			\scaleto{%
				\scalerel*[\widthof{\ensuremath{#1}}]{\kern-.6pt\bigwedge\kern-.6pt}%
				{\rule[-\textheight/2]{1ex}{\textheight}}
			}{\textheight}%
		}{0.5ex}}%
	\stackon[1pt]{#1}{\tmpbox}%
}

\DeclareDocumentCommand\issue{g}{\todo[size=\footnotesize,color = green!40]{Issue\IfNoValueF{#1}{: #1}}}
\DeclareDocumentCommand\tobedone{g}{\todo[size=\footnotesize,color = yellow!50]{To do\IfNoValueF{#1}{: #1}}}
\DeclareDocumentCommand\notationissue{g}{\todo[size=\footnotesize,color = red!30]{Notation?\IfNoValueF{#1}{: #1}}}
\DeclareDocumentCommand\doubt{g}{\todo[size=\footnotesize,color = blue!10]{Doubt\IfNoValueF{#1}{: #1}}}
\DeclareDocumentCommand\observation{g}{\todo[size=\footnotesize,color = orange!10]{Observation\IfNoValueF{#1}{: #1}}}

\usepackage[all,cmtip]{xy}

\usetikzlibrary{calc}
\usetikzlibrary{matrix,arrows,arrows.meta}
\usepackage{tikz-cd}

\usepackage{calligra} 

\newtagform{fn}{(}{)\footnotemark}

    {\endMakeFramed}

    {\endMakeFramed}

\makeatletter
\g@addto@macro\bfseries{\boldmath}
\makeatother

\theoremstyle{definition}
\newtheorem{thm}{Theorem}[section]
\newtheorem*{thm*}{Theorem}
\newtheorem{prop}[thm]{Proposition}

\newtheorem{cor}[thm]{Corollary}
\newtheorem{constr}[thm]{Construction}

\theoremstyle{definition}
\newtheorem{defn}[thm]{Definition}

\theoremstyle{remark}
\newtheorem{rem}[thm]{Remark}

\newcommand\cu{\mathcal {U}}

\newcommand\gh{\mathscr{GH}}
\newcommand\gsh{\mathscr{SH}_G}

\newcommand\m{\mathbf{MU}}
\newcommand\kk{\mathbf{KU}}

\newcommand{\Ab}{\ensuremath{\mathsf{Ab}}}

\mathchardef\mhyphen="2D

\title[Completion and Local cohomology for Compact Lie Groups]{The completion and local cohomology theorems for complex cobordism for all compact Lie groups}
\author{Marco \textsc{La Vecchia}} 

\address{
\newline Marco \textsc{La Vecchia};
University of Warwick;
Mathematics Institute;
Coventry, United Kingdom
\textit{Email adress:} \href{mailto:marco.la-vecchia@warwick.ac.uk}{\sf marco.la-vecchia@warwick.ac.uk}\newline
}

\begin{document}
\tikzcdset{arrow style=tikz, diagrams={>=to}}

\begin{abstract}  We generalize the completion theorem for equivariant $\m$-module spectra for finite groups or finite extensions of a torus to compact Lie groups using the splitting of global functors proved by Schwede. 
This proves a conjecture of Greenlees and May.
\end{abstract}
\maketitle

\tableofcontents

\section{Overview}

\subsection{Introduction}

A completion theorem establishes a close relationship between equivariant cohomology theory and its non-equivariant counterpart. It takes various forms, but in favourable cases it states
\begin{equation*}
(E^*_G)^{\wedge}_{J_G} \cong E^*(BG)
\end{equation*}
where $E$ is a $G$-spectrum, $E^*_G$ is the associated equivariant cohomology theory and $J_G$ is the augmentation ideal (\Cref{augmentationIdeal}).

The first such theorem is the Atiyah-Segal Completion Theorem for complex K-theory \cite{AtyahSegal}. This especially favourable because the coefficient ring $\kk_G^*=R(G)[v, v^{-1}]$ is well understood and in particular it is Noetherian, and so in this case we can view the theorem as the calculation of the cohomology of the classifying space.  The good behaviour for all groups permits one to make good use of naturality in the group, and indeed \cite{AtyahSegal} uses this to give a proof  for all compact Lie groups $G$. The result raised the question of what other theories enjoy a completion theorem, and the case of equivariant complex cobordism was considered soon afterwards, with Löffler giving a proof in the abelian case \cite{Loffler:1974uz}.  The fact that the coefficient ring $\m_G^*$ is not known explicitly means that this cannot be viewed as a computation of the cohomology of the classifying space. The fact that the coefficient ring is unknown and not Noetherian was an obstacle to extensions.

 Despite the algebraic complexity of the coefficients, Segal made the remarkable conjecture that stable cohomotopy should satisfy the completion theorem,  and this was proved by Carlsson \cite{CarlssonSegal}, building on important earlier work by others (see, e.g., \cite{lin_1980}, \cite{lin_davis_mahowald_adams_1980}, \cite{Ravenel81thesegal}, \cite{laitinen} \cite{segal_stretch_1981},\cite{stretch_1981}, \cite{CARLSSON198383}, \cite{ADAMS1985435}). In this case the conclusion can only be viewed as a calculation of the cohomotopy of classifying spaces in degrees 0 and below, but the structural content in positive degrees is equally striking. In the course of understanding this, there was a focus on understanding completion in various ways. From a homotopical point of view this led to the connection between completion and local cohomology and the definition of local homology \cite{GREENLEES1992438}. It also led to a new proof of the Atiyah-Segal Completion Theorem and also its counterpart in homology \cite{GREENLEES1993295}. This in turn reopened the question of the completion theorem and local cohomology theorem for $\m$ but now with the challenges shifted from the formal behaviour to the algebraic behaviour: the formal structure of the proof of the local cohomology theorem for $\kk$ applies precisely for $\m$, but the difficulty is that since $\m_G^*$ is not Noetherian, it is not clear that the relevant ideals are finitely generated. Accordingly, Greenlees and May (\cite{GreenleesCompletion}) isolated the formal argument and observed that if one could find a sufficiently large finitely generated ideal (\Cref{SUFFLARGE}), the local cohomology and completion theorems would hold for $\m$. They went on to codify and use the structure of $\m$ as a global spectrum to define and apply multiplicative norm maps and hence construct sufficiently large finitely generated ideals in the case when the identity component of $G$ is a torus. This led to the proof of  local cohomology and completion theorems for $\m$ for these groups.

Much more recently, Schwede has studied global spectra more systematically \cite{schwede_2018}, and in particular used the global structure of $\m$ in a more sophisticated way to show that tautological unitary Euler classes are regular and give rise to various splittings \cite{schwede2020splittings}. 
The contribution of the present paper is to show that every compact Lie group $G$ has a sufficiently large finitely generated ideal, and in fact if $G$ embeds in $U(n)$ it gives generators in terms of Euler classes of the tautological representations of $U(k)$ for $k\leq n$. This is a consequence of Schwede's results as we will see in \Cref{ourthm}.

\begin{thm}
For any compact Lie group $G$ with a faithful representation in $U(n)$ and any $G$-space $X$ there are spectral sequences
$$H_{J_G}^*(\m^G_*(X))\Rightarrow \m^G_*(EG_+\wedge X)$$
$$H^{J_G}_*(\m_G^*(X))\Rightarrow \m_G^*(EG_+\wedge X), $$
where the local homology and cohomology are calculated using a subideal of $J_G$ with $n$-generators and therefore concentrated in degrees $\leq n$. 
\end{thm}

\begin{rem} \label{homologycompletion}
For a Noetherian ring $R$ and an ideal $I$ the local homology groups calculate the left derived functors of completion \cite{GREENLEES1992438}
$$H_*^I(R;M)=L^I_*M,$$
and the local cohomology groups 
calculate the right derived functors of $I$-power torsion \cite{Grothendieck}
$$H^*_I(R;M)=R^*\Gamma_I(M). $$
This is not known when $R=\m^G_*$ and $J=J_G$, but the local homology and cohomology with respect to $n$ specified elements are much  more practical to  compute.
\end{rem}

\subsection{Organization} We start with a preliminary section where we introduce the notation and some basic facts of equivariant and global orthogonal spectra. In \Cref{GreenleesCompletion} we review the classical statement. This section is only needed to recall basic constructions and state the main theorem that we will prove in \Cref{ourthm}. In \Cref{SchwedeSplitting} we review Schwede's splitting (\cite[Theorem 1.4, p. 5]{schwede2020splittings}) and his corollary that ensures the regularity of certain Euler classes \cite[Corollary 3.2, p.10]{schwede2020splittings}.
Finally, in \Cref{ourthm} we prove that the augmentation ideal $J_{U(n)}$ is \textit{sufficiently large}. This will implies the completion theorem \cite[Theorem 1.3, p. 514]{GreenleesCompletion} for $U(n)$ and for any compact Lie group $G$.

\subsection{Acknowledgment}I am extremely grateful to my supervisor Prof. Dr. John Greenlees, for introducing me to the problem and for his constant support, advice and help during the project. I also would like to thank Prof. Dr. Stefan Schwede for reading a preliminary version of the article and suggesting an improvement on the proof of the main result.

\section{Notations, Conventions and Facts} \label{prel}
\subsection{Spaces}
By a \textit{space} we mean a \textit{compactly generated space} as introduced in \cite{McCord1969ClassifyingSA}. We will denote by $\mathscr{Top}$ (respectively $\mathscr{Top}_*$)  the category of spaces (pointed spaces) with continuous maps (based maps). 
For a compact Lie group $G$, $\mathscr{Top}_G$ denotes the category of $G$-spaces and $G$-equivariant maps. (Equivariant) mapping spaces and (based) homotopy classes of (based) maps are defined as usual and are denoted by $\textrm{map}(\_, \_)$ and $[\_, \_]$ respectively.
\subsubsection{Universal spaces}\label{UniversalSpaces} A \textit{family} of subgroups of a group $G$ is a collection of subgroups closed under conjugation and taking subgroups.
When $G$ is a compact Lie group, a universal $G$-space for a family $\mathcal{F}$ of closed subgroups is a cofibrant $G$-space $E\mathcal{F}$ such that:
\begin{itemize}
	\item all isotropy groups of $E\mathcal{F}$ belong to the family $\mathcal{F}$, 
	\item for every $H \in \mathcal{F}$ the space $E\mathcal{F}^H$ is weakly contractible.
\end{itemize}
We denote by $\tilde{E}\mathcal{F}$  the reduced mapping cone of the collapse map $E\mathcal{F}_+ \rightarrow S^0$ which sends $E\mathcal{F}$ to the non basepoint of $S^0$.
Any two such universal $G$-spaces are $G$-homotopy equivalent, hence we will refer to $E\mathcal{F}$  as the universal space for the family $\mathcal{F}$.
Note that if the family $\mathcal{F}= \{e\}$ then $E\mathcal{F}= EG$.

\subsection{Algebra} \label{Loccohomology} For a graded commutative ring $A$ and a finitely generated ideal $I= (a_1, \dots, a_n)$ 
of $A$ we let $K_{\infty}^{\bullet}(I)$ be the graded cochain complex  
\begin{equation*}
\bigotimes_{i=1, \dots, n}(A\rightarrow A[1/a_i]),
\end{equation*}
where $A$ and $A[1/a_i]$ sit in homological degrees $0$ and $1$ respectively. If $N$ is a graded $A$-module the local cohomology groups are defined as
\begin{equation*}
H^{s,t}_I(A; N)= H^{s,t}(K_{\infty}^{\bullet}(I)\otimes N).
\end{equation*}
When $A$ is Noetherian the functor $H^*_I(A; \_)$ calculates the right derived functors of the torsion functor \begin{equation*}
\Gamma_I(N)= \{n \in N \: \textrm{such that} \: I^kn=0 \: \textrm{for some} \: k\}.
\end{equation*}

The main references for the theory of local cohomology are \cite{GrothendieckFR} and \cite{Grothendieck}. 

Dually we let the local homology groups be \begin{equation*}
H^I_{s,t}(A; N)= H_{s,t}(\textrm{Hom}(\tilde{K}_{\infty}^{\bullet}(I), N),
\end{equation*}
where $\tilde{K}_{\infty}^{\bullet}(I)$ is a $A$-free chain complex quasi-isomorphic to $K_{\infty}^{\bullet}(I)$ (\cite{GREENLEES1992438}).
As already mentioned in \Cref{homologycompletion}, when $A$ is Noetherian and $N$ is free or finitely generated then the functor $H^I_k(A; \_)$ calculates the left derived functors of the $I$-adic completion functor. In particular, under these assumptions
\begin{equation*}
H^I_k(A; N) \cong
\begin{cases}
 N^{\wedge}_I \quad \textrm{if} \: k=0; \\
0 \quad \textrm{otherwise}.
\end{cases}
\end{equation*}

\subsection{Spectra} A spectrum will be an orthogonal spectrum and we will denote by $\textrm{Sp}$ the category of orthogonal spectra as defined in \cite[Definition 3.1.3, p.230]{schwede_2018}. The category $\textrm{Sp}$ is closed symmetric monoidal with respect to the smash product $\_ \wedge \_$. We will denote by $\textrm{map(\_, \_)}$ the right adjoint of the smash product. $\textrm{Sp}$ is also tensored over $\mathscr{T}_*$, i.e., for every based space $A$ and every spectrum $X$ a mapping spectrum between these two is defined. We will use the notation $\textrm{map(\_, \_)}$ also in this case. The definitions in the equivariant case are similar.

\subsection{Global and Equivariant Stable Homotopy Categories}
We let $\gh$ and $\gsh$ denote respectively the \textit{global stable homotopy category} and the \textit{$G$-equivariant stable homotopy category} for a compact Lie group $G$. The first can be realized as a localization of the category of orthogonal spectra at the class of \textit{global equivalences} \cite[Definition 4.1.3, p. 352]{schwede_2018} as constructed in \cite[Theorem 4.3.18, p. 400]{schwede_2018}. The second one as a localization at the class of $\pi_*$-\textit{isomorphisms} of the category of $G$-orthogonal spectra as constructed in \cite[Theorem 4.2, p. 47]{MandellMay}.
Hence, we have 
\begin{align*}
\gsh & \cong \textrm{Sp}_G[(\pi_*\textrm{-isos})^{-1}] \\
\gh \cong  \textrm{Sp}&[(\textrm{global equivalences})^{-1}].
\end{align*}

 Both global equivalences and $\pi_*$-isomorphisms are weak equivalences of stable model structures, see \cite[Theorem 4.3.17, p. 398]{schwede_2018} for the global case and \cite[Theorem 4.2, p. 47]{MandellMay} for the equivariant one. This implies that both categories $\gh$ and $\gsh$ come with a prefered structure of triangulated categories and we denote by $\Sigma$ the \textit{shift functor} in both cases.  
The derived smash product of $\textrm{Sp}$ (resp. $\textrm{Sp}_G$) endows the category $\gh$ (resp. $\gsh$) with a closed symmetric monoidal structure.

Homotopy groups are defined for equivariant spectra and for global spectra as usual  \cite[p. 232]{schwede_2018}. In both cases for a fixed $X$ and $k$, the system of homotopy groups $\{ \pi_k^H(X)\}_{H \subset G}$ has a lot of additional structure. For $G$-spectra $\{ \pi_k^H(X)\}_{H \subset G}$ is a \textit{Mackey functor} (\cite[Definition 3.4.15, p. 319]{schwede_2018}), for global spectra it is a \textit{global functor} (\cite[Definition 4.2.2, p. 369]{schwede_2018}).

If we fix a compact Lie group G and a $G$-spectrum $X$, the functor 
\begin{equation*}
\pi_k^G \colon \gsh \rightarrow \Ab
\end{equation*}
is represented by the pair $(\Sigma^{k}\mathbb{S}, \textrm{id})$.
Hence, we have a natural isomorphism

 \begin{equation}\label{representedFunctor}
 \gsh(\Sigma^k\mathbb{S}, X) \cong \pi_k^G(X).
 \end{equation}

 A similar statement holds in the global setting, refer to \cite[Theorem 4.4.3, p. 412]{schwede_2018}.
 Finally, for a $G$-spectrum $X$ we adopt the usual convention 

 \begin{equation*}
 X^G_*= \pi_*^G(X), \qquad  X^*_G= X^G_{-*}.
 \end{equation*}

\subsection{Global Complex Cobordism} \label{MU} Since our work regards the complex bordism ring, we recollect here some facts about this theory.
We will write $\m$ for the global Thom ring spectrum defined in \cite[ Example 6.1.53]{schwede_2018} which represents the complex bordism spectrum nonequivariantly, and for a fixed compact Lie group $G$, it is a model for the homotopical equivariant bordism introduced by tom Dieck in \cite{TOMDIECK1970345}. The global theory $\m$ has the structure of an \textit{ultracommutative ring spectrum} in the sense of \cite[Definition 5.1.1, p. 463]{schwede_2018}. 
For every compact Lie group $G$ and for every unitary representation $V$, a \textit{Thom class} $\sigma_G(V) \in \m^{2n}_{G}(S^V)$ is defined, where $n=\textrm{dim}_{\mathbb{C}}(V)$. The Euler class $e_{G}(V)$ is by definition the image of the Thom class along the fixed point inclusion $S^0 \rightarrow S^V$. If $V$ has non trivial $G$-fixed-points, then the previous inclusion is $G$-null-homotopic hence, $e_{G}(V)=0$ if $V^G \neq \{0\}$. 
On the other hand tom Dieck showed that if $V^G = \{0\}$ then the Euler class $e_G(V)$ is a non zero element in $\m^{2n}_G$ \cite[Corollary 3.2, p. 352]{TOMDIECK1970345}.

The theory $\m$ has an equivariant Thom isomorphism for every unitary representation $V$ 
\begin{equation}
\m^0_G(S^{k+V}) \cong \m^{-k-2n}_G
\end{equation}
where $n= \textrm{dim}_{\mathbb{C}}(V)$. This ismorphism takes the multiplication by the class $a_n$ defined in \Cref{SchwedeSplitting} to multiplication by the Euler class of the representation $V$.

\section{Classical Statement}\label{GreenleesCompletion}

We recall some basic constructions that can be found in \cite{GreenleesCompletion}.

Let $R$ be an orthogonal $G$-ring spectrum, i.e., a monoid in the monoidal category of orthgonal $G$-spectra with respect to the smash product (\cite[Section 3.5, p. 333]{schwede_2018}, \cite[Chapter II, section 3, p. 33]{MandellMay}). 
We could assume $R$ to be a $G$-ring spectrum in a weaker sense but, since we are mainly interested in studying  the global Thom spectrum $\m$ which has the structure of an \textit{ultracommutative ring spectrum} \cite[Definition 5.1.1, p. 463]{schwede_2018} we use this stronger notion to simplify the exposition.

\begin{constr}\label{LocHom} Let $R$ be as above. By \Cref{representedFunctor}, every element of $R^G_n$ specifies by adjunction a morphism $\alpha \colon \mathbb{S} \rightarrow \Sigma^{-n}R$ in $\gsh$. 
We let 

\begin{equation*}
\tilde{\alpha} \colon R \rightarrow \Sigma^{-n}R. 
\end{equation*}

\hspace{-17pt} be the following composition

\begin{equation*}
R \xrightarrow{\alpha\wedge R} \Sigma^{-n}R \wedge R\cong \Sigma^{-n}(R\wedge R) \xrightarrow{\Sigma^{-n}\mu}  \Sigma^{-n}R,
\end{equation*}

\hspace{-17pt} where $\mu \colon R\wedge R \rightarrow R$ is the multiplication of $R$.
This defines a morphism in $\gsh$ and we let

\begin{equation}
R[1/ \alpha] := \textrm{telescope}(R\xrightarrow{\tilde{\alpha}} \Sigma^{-n}R \xrightarrow{\Sigma^{-n}\tilde{\alpha}} \Sigma^{-2n}R\rightarrow \cdots ),
\end{equation}

\hspace{-17pt} be the \textit{mapping telescope} of the iterates of $\tilde{\alpha}$.
\end{constr}

\begin{rem}The \textit{mapping telescope} in \Cref{LocHom} models the sequential homotopy colimit in $\gsh$. 
For a discussion of homotopy colimits in $\gsh$ refer to \cite[Appendix C, p. 160]{NikSch}.
\end{rem}

\begin{defn}\label{Koszul} Let $R$ and $\alpha$ be as above and let $I = (\alpha_1, \dots , \alpha_n)$ be an ideal in $R^G_*$. We define \begin{align*}
&K_{\infty}(\alpha):=\textrm{fib}(R \rightarrow R[1/ \alpha]), \\
K_{\infty}(&I):= K(\alpha_1) \wedge_{R} \dots \wedge_{R}K(\alpha_n).
\end{align*}
\end{defn}

\begin{defn}\label{GAMMA} Let $M$ be an $R$-module and $I \subset R^G_*$ be a finitely generated ideal. Then we define
\begin{equation*}
\Gamma_{I}(M)= K_{\infty}(I)\wedge_R M
\end{equation*}
and
\begin{equation*}
(M)^{\wedge}_{I}= \textrm{map}_{R}(K_{\infty}(I), M).
\end{equation*}
\end{defn}

\begin{rem}\label{spectralsequence} There is a spectral sequence of local cohomology
\begin{equation*}
H^*_I(R^G_*; M^G_*) \Rightarrow \Gamma_I(M)^G_*,
\end{equation*}
and there is a spectral sequence of local homology 
\begin{equation*}
H^I_*(R^G_*; M^G_*) \Rightarrow ((M)^{\wedge}_{I})^*_G.
\end{equation*}
Note that when $M=R$ we obtain the spectral sequence that computes $K_{\infty}(I)^G_*$.
\end{rem}

{}
We now turn our attention to $\m$ (\Cref{MU}).

\begin{defn}\label{SUFFLARGE}(\cite[Definition 2.4, p. 517]{GreenleesCompletion}) An ideal $I \subset \m^*_G$ is sufficiently large at $H$ if there exists a nonzero complex representation $V$ of $H$ such that $V^H=\{0\}$ and the Euler class $e_H(V) \in \m^{2n}_{H}$ is in the radical $\sqrt{res^G_H(I)}$, where $n=\textrm{dim}_{\mathbb{C}}(V)$. The ideal $I$ is sufficiently large if it is sufficiently large at all closed subgroup $H\neq 1$ of $G$. 
\end{defn}

\begin{defn}\label{augmentationIdeal} Let $G$ be a compact Lie group and $R$ be an orthogonal $G$-spectrum. The augmentation ideal of $R$ at $G$ is the kernel of \begin{equation*}
\textrm{res}^G_1 \colon R^G_* \rightarrow R_*.
\end{equation*}
\end{defn}

\begin{constr}\label{KAPPA}By construction of $R[1/\alpha]$  if $\alpha \in J_G$ then

 \begin{equation*}
\textrm{res}^G_{1} R[1/\alpha] \simeq 0. 
\end{equation*} 

\hspace{0pt} Hence,  applying the restriction to the fibre sequence 
\[\begin{tikzcd}
	{\Gamma_{\alpha} R} & R & {R[1/\alpha]}
	\arrow[from=1-1, to=1-2]
	\arrow[from=1-2, to=1-3]
\end{tikzcd}\]

\hspace{-17pt} we obtain a fibre sequence in which the third term is contractible. This implies, by the long exact sequence in homotopy groups induced by a fibre sequence of spectra, that the canonical map \begin{equation*}
\textrm{res}^G_1 (\Gamma_{\alpha} R)=\Gamma_{\textrm{res}^G_1(\alpha)} \textrm{res}^G_1 R \xrightarrow{\simeq} \textrm{res}^G_1 R
\end{equation*}
is an equivalence.
The same argument applies for an ideal $I \subset R^G_*$, obtaining an equivalence
\begin{equation*}
\textrm{res}^G_1 (\Gamma_{I} R)=\Gamma_{\textrm{res}^G_1(I)} \textrm{res}^G_1 R \xrightarrow{\simeq} \textrm{res}^G_1 R.
\end{equation*}
Smashing with the universal $G$-space $EG_+$ the above morphism, we obtain an equivalence
 \begin{equation*}
EG_+ \wedge \Gamma_{I} R \rightarrow EG_+\wedge R
\end{equation*}
in $\gsh$. 
Inverting this equivalence and composing with the collapse map $EG_+ \wedge \Gamma_{I} R \xrightarrow{\textrm{coll} \wedge \Gamma_{I} R} S^0 \wedge \Gamma_{I} R \cong \Gamma_{I} R$ we obtain a zig-zag
\[\begin{tikzcd}
	{EG_+\wedge R} & {EG_+\wedge \Gamma_{I} R} & {\Gamma_{I} R}
	\arrow["\simeq"', from=1-2, to=1-1]
	\arrow[from=1-2, to=1-3]
	\arrow["\kappa", curve={height=-24pt}, from=1-1, to=1-3],
\end{tikzcd}\]
which defines a morphism of $R$-modules in $\gsh$. 
\end{constr}

\begin{thm}(\cite[Theorem 2.5, p.518]{GreenleesCompletion})\label{GreenLoc}
Let $G$ be a compact Lie group. Then, for any sufficiently large finitely generated ideal $I \subset J_G$
\begin{equation*}
\kappa \colon EG_+ \wedge \m \rightarrow \Gamma_I\m
\end{equation*}{}
is an equivalence in $\gsh$. Therefore, 
\begin{equation*}
EG_+ \wedge M \rightarrow \Gamma_{I}(M) \quad \textrm{and} \quad (M)^{\wedge}_{I} \rightarrow \textrm{map}(EG_+, M) 
\end{equation*}
are equivalences for any $\m$-module $M$.
\end{thm}

\begin{proof} Here, we only give a sketch of the argument following the main reference. The point is that if $I \subset \m^*_G$ is sufficiently large then $\textrm{res}^G_HI \subset \m^*_H$ is also sufficiently large. Moreover, since every descending sequence of compact Lie groups stabilizes, we can use induction and assume that the theorem holds for any proper closed subgroup of $G$. Passing to the cofibre of the map $\kappa$ it is enough to show that \begin{equation*}
\pi_*^G(\tilde{E}G \wedge \Gamma_I\m)=0.
\end{equation*}
We then let $\mathcal{P}$ be the family of proper subgroup of $G$ and let $E\mathcal{P}$ be the universal space associated to $\mathcal{P}$. Since \begin{equation*}
\tilde{E}\mathcal{P} \wedge S^0 \rightarrow \tilde{E}\mathcal{P} \wedge \tilde{E}G
\end{equation*} 
is an equivalence, it suffices to show that $\tilde{E}\mathcal{P}\wedge K(I)$ is contractible. 
Let $\cu$ be a complete complex $G$-universe and define $\cu^{\perp}$ to be the orthogonal complement of the $G$-fixed points $\cu^G$ in $\cu$. Then, \begin{equation*}
\textrm{colim}_{V \in \cu^{\perp}} S^V
\end{equation*}
is a model for $\tilde{E}\mathcal{P}$. 
We can compute then 
\begin{align*}
&\pi_*^G  (\tilde{E}\mathcal{P} \wedge \Gamma_I\m) = \pi_*^G((\textrm{colim}_{V \in \cu^{\perp}} S^V)\wedge \Gamma_I\m) \cong \\
   \textrm{col}&\textrm{im}_{V \in \cu^{\perp}} \pi_*^G( S^V\wedge \Gamma_I\m )  \cong \textrm{colim}_{V \in \cu^{\perp}} \pi^G_{*-|V|}(\Gamma_I\m) \cong  \\ 
   & \qquad \quad \pi_*^G(\Gamma_I\m) [\{e_G(V)^{-1}\}_{V \in \cu^{\perp}}] .
\end{align*}

Localising the spectral sequence in \Cref{spectralsequence} \begin{equation*}
H_I^*(\m^G_*) \Rightarrow \pi_*^G(\Gamma_I\m)
\end{equation*}
away from the Euler classes we obtain another spectral sequence 
\begin{equation*}
H_I^*(\m^G_*)[\{e_G(V)^{-1}\}_{V \in \cu^{\perp}}] \Rightarrow \pi_*^G(\Gamma_I\m )[\{e_G(V)^{-1}\}_{V \in \cu^{\perp}}].
\end{equation*}
Since local cohomology of a ring at an ideal become zero when localized by inverting an element in that ideal, we obtain that the $E^2$-term of the spectral sequence is zero being $I$ sufficently large. This proves the claim.
\end{proof}

\begin{rem} As stated in \cite[Theorem 2.5, p.518]{GreenleesCompletion}, the previous theorem holds for any equivariant ring spectrum $R$ with Thom isomorphism. Examples of these theories are the homotopical equivariant bordism $\m$ and the equivariant $K$-theory. For the purpose of this paper we should require that $R$ is an ultracommutative ring spectrum with Thom isomorphism.
\end{rem}

The paper \cite{GreenleesCompletion} proceeds by constructing a sufficiently large subideal of the augmentation ideal $J_G$ whenever $G$ is a finite group or a finite extension of a torus. Hence, if we find a sufficiently large subideal of the augmentation ideal $J_G$ for any compact Lie group, we obtain the Completion \Cref{GreenLoc} for any compact Lie group. We use Schwede's Splitting (\ref{SchwedeSplit}) to prove that $J_{U(n)}$ is generated by "Euler classes". Thanks to this we prove that $J_{U(n)}$ is sufficiently large (\Cref{OURTHM}). Using the fact that any compact Lie group embeds into a unitary group $U(N)$ for $N$ sufficiently large we conclude that $\textrm{res}^{U(N)}_G J_{U(N)}$ is a sufficiently large subideal of $J_G$ for any compact Lie group $G$.

\section{Schwede's Splitting} \label{SchwedeSplitting}

We recall that a \textit{global functor} $F$ associates to every compact Lie group $G$ an abelian group $F(G)$ and this association is controvariantly functorial with respect to continous group homomorphisms. Moreover, for every closed subgroup inclusion $H < G$  a transfer map $\textrm{tr}^G_H \colon M(H) \rightarrow M(G)$ is defined. This data needs to satisfy some relations that can be found in \cite[p. 373]{schwede_2018}.
In \cite[Theorem 1.4, p. 5]{schwede2020splittings}, Schwede proves that for any global functor $F$ the restriction homomorphism 
\begin{equation*}
\textrm{res}^{U(n)}_{U(n-1)} \colon F(U(n)) \rightarrow F(U(n-1))
\end{equation*}
is a split epimorphism. He then deduces a splitting of global functors when evaluated on the unitary group $U(n)$.
Explicitly, the splitting take the form \begin{equation}\label{SchwedeSplit}
F(U(n)) \cong F(e) \bigoplus \bigoplus_{k=1, \dots , n} \textrm{Ker}(\textrm{res}^{U(k)}_{U(k-1)} \colon F(U(k)) \rightarrow  F(U(k-1))).
\end{equation}

\begin{rem}\label{forget}There is a forgetful functor $U\colon \gh \rightarrow \gsh$ which is strong symmetric monoidal and exact. Hence, the global splitting (\ref{SchwedeSplit}) at the unitary group translates in a splitting in $\gsh$.
\end{rem}

 The most important application of the splitting for us is when the global functor comes from the homotopy groups of an orthogonal spectrum. In fact, for every global stable homotopy type $X$, i.e., an object in $\gh$, we have a global functor \begin{equation*}
\underline{\pi}_*(X)(G)= \pi_*^G(X).
\end{equation*}
The splitting then tells us that for every $k \le n$ the group $\pi_*^{U(k)}(X)$ is a natural summand of $\pi_*^{U(n)}(X)$. 
In this case a more explicit description of the right hand side of the splitting is available. 
In fact, let $\nu_n$ be the tautological representation of $U(n)$ and let
\begin{equation*} 
a_n \in \pi_0^{U(n)}(\Sigma^{\infty} S^{\nu_n}) 
\end{equation*} be the Euler class of $\nu_n$, i.e., the element represented by the inclusion $S^0 \rightarrow S^{\nu_n}$. Then the following short sequence

\[\begin{tikzcd}\label{EULERes}
	0 & {\pi_{* + \nu_n}^{U(n)}(X)} & {\pi_{*}^{U(n)}(X)} & {\pi_{*}^{U(n-1)}(X)} & 0
	\arrow["{\textrm{res}^{U(n)}_{U(n-1)}}", from=1-3, to=1-4]
	\arrow[from=1-4, to=1-5]
	\arrow[from=1-1, to=1-2]
	\arrow["{a_n}", from=1-2, to=1-3]
\end{tikzcd}\]
is exact \cite[Corollary 3.1, p. 10]{schwede2020splittings}.

When $X=\m$ the equivariant Thom isomorphism identifies the short exact sequence (\ref{EULERes}) in the following short exact sequence

\[\begin{tikzcd}
	0 & {\m^{*-2n}_{U(n)}} & {\m^*_{U(n)}} & {\m^*_{U(n-1)}} & 0
	\arrow[from=1-1, to=1-2]
	\arrow["{e_{U(n)}(\nu_n)}", from=1-2, to=1-3]
	\arrow["{\textrm{res}^{U(n)}_{U(n-1)}}", from=1-3, to=1-4]
	\arrow[from=1-4, to=1-5].
\end{tikzcd}\]

Moreover, we have the following corollary.

\begin{cor}\label{AUGMENTATIONN} Let $J_{U(n)}$ be the augmentation ideal of $\m^{U(n)}_*$, i.e., the kernel of the $\textrm{res}^{U(n)}_{1}\colon  \m^{U(n)}_* \rightarrow \m^{1}_*$, and let $s^{U(n)}_{U(k)}$ be a section of $\textrm{res}^{U(n)}_{U(k)}$ (see \cite[Construction 1.3, p.4]{schwede2020splittings}). Then\begin{equation*}
J_{U(n)}= (s^{U(n)}_{U(k)}(e_{U(k)}(\nu_k))  \quad|\quad  k=1, \dots, n)
\end{equation*}
and in particular, 
\begin{equation*}
\textrm{res}^{U(n)}_{U(k)}J_{U(n)}= J_{U(k)}
\end{equation*}
for all $k\leq n$.
\end{cor}
\begin{proof}This is just the combination of the splitting \ref{SchwedeSplit} and the short exact sequence \ref{EULERes}.
\end{proof}

\section{The main Result}\label{ourthm}

\begin{prop}\label{OURTHM}The augmentation ideal of $\m$ at the unitary group $U(n)$ is sufficiently large.
\end{prop}

\begin{proof}Let $H$ be a closed subgroup of $U(n)$. Consider $$V=\textrm{res}^{U(n)}_{H} \nu_n - (\textrm{res}^{U(n)}_{H} \nu_n)^H$$ and let $k=\textrm{dim}_{\mathbb{C}}(V)$. Note that $k=0$ if and only if $H=1$. 
We claim that the Euler class $e_H(V)$ is in the $\textrm{res}^{U(n)}_H J_{U(n)}$.
If $k=n$ then $V=\textrm{res}^{U(n)}_{H} \nu_n$ and \begin{equation*}
e_H(V)= \textrm{res}^{U(n)}_{H} (e_{U(n)} \nu_n) \in \textrm{res}^{U(n)}_H J_{U(n)}
\end{equation*}
and the claim holds.

\hspace{-16pt} If otherwise $k\neq 0$, we choose an  orthonormal basis $(x_1, \dots, x_{n-k})$ of $(\textrm{res}^{U(n)}_{H} \nu_n)^H$ and a unitary matrix $g \in U(n)$ that sends the canonical basis of $\mathbb{C}^n$ to any other orthonormal basis that has as last $n-k$ vectors $(x_1, \dots, x_{n-k})$.
Then for any $h\in H$

\begin{equation*}
h^g=\left(\begin{array}{c|c}
  \raisebox{-20pt}{\qquad $\tilde{h}$ \qquad} & \raisebox{-20pt}{$0$} \\ 
& \, \\ \hline
 \raisebox{-10pt}{$0$} & \raisebox{-10pt}{\,\,\,$\textrm{Id}_{n-k}$ } \\ 
& \, \\
  \end{array} \right).
\end{equation*} \\

\hspace{-16pt} This implies that $V$ is conjugate to the $H^g$-representation $\textrm{res}^{U(k)}_{H^g}(\nu_k)$.
Letting  $g_{\star} \colon \m^*_{H^g} \rightarrow \m^*_{H}$ be the conjugation action  (see relations in \cite[Definition 3.4.15, p. 319]{schwede_2018}), we pass to Euler classes obtaining the relation \begin{equation*}
e_H(V)= g_{\star} (e_{H^g}(\textrm{res}^{U(k)}_{H^g}(\nu_k))).
\end{equation*}
We then compute the right hand side of the last equation
\begin{align*}
g_{\star} (e_{H^g}(\textrm{res}^{U(k)}_{H^g}(\nu_k)))= g_{\star} (\textrm{res}^{U(k)}_{H^g}(e_{U(k)}(\nu_k)))= \\
g_{\star} (\textrm{res}^{U(n)}_{H^g}( s^{U(n)}_{U(k)}(e_{U(k)}(\nu_k))))=\textrm{res}^{U(n)}_H ( s^{U(n)}_{U(k)}(e_{U(k)}(\nu_k)))
\end{align*}
where, in the second equality we have used the chosen section $s^{U(n)}_{U(k)}$ (\Cref{AUGMENTATIONN}) and in the last one the formula \begin{equation*}
g_{\star} \circ \textrm{res}^{U(n)}_{H^g}= \textrm{res}^{U(n)}_H
\end{equation*} (again, see relations in \cite[Definition 3.4.15, p. 319]{schwede_2018}).
By \Cref{AUGMENTATIONN} it is clear that $\textrm{res}^{U(n)}_H ( s^{U(n)}_{U(k)}(e_{U(k)}(\nu_k))) \in \textrm{res}^{U(n)}_H J_{U(n)}$, hence we have proved the claim. 
Since by contruction $e_H(V)$ is non zero we conclude that $J_{U(n)}$ is sufficiently large at $H$.

\end{proof}

We now let $G$ be any compact Lie group. Since every compact Lie group has a faithful representation, $G$ is isomorphic to a closed subgroup of $U(n)$ where $n$ is the dimension of a chosen faithful representation of $G$. Then we have the following corollary.

\begin{cor}
The ideal $\textrm{res}^{U(n)}_G J_{U(n)} \subset J_G$ is a sufficiently large finitely generated ideal. In particular, the completion \Cref{GreenLoc} holds for any compact Lie group $G$ if we choose $I= \textrm{res}^{U(n)}_G J_{U(n)}$.
\end{cor}
\begin{proof}
By transitivity of restrictions the ideal $\textrm{res}^{U(n)}_G J_{U(n)}$ is sufficiently large and is contained in $J_G$. The completion theorem then holds by the argument above.
\end{proof}

\begin{rem} The subideal $J= \textrm{res}^{U(n)}_G J_{U(n)}$ of $J_G$ is not special. Indeed, if $I$ is any other finitely generated subideal of $J_G$ containing $J$, then 
\begin{equation*}
\Gamma_{I}\m \simeq \Gamma_{J} \m,
\end{equation*}
 and

\begin{equation*}
\m^{\wedge}_{I} \simeq \m^{\wedge}_{J}.
\end{equation*}
In fact, \Cref{GreenLoc} implies that the $\m$-modules $K_{\infty}(I)$ and $\Gamma_{I}M$ are independent of the choice of $I$.
\end{rem}

\bibliography{thebib}
\bibliographystyle{alpha}

\end{document}